\documentclass[11pt,a4paper]{amsart}

\usepackage{epsfig}
\usepackage{amssymb,amsmath,amsfonts}
\usepackage{amsthm,amstext,array}
\usepackage{latexsym}
\usepackage{bbm}
\usepackage{graphicx}
\usepackage{mathtools}
\usepackage{enumerate}
\usepackage{MnSymbol}
\usepackage{caption}
\usepackage{subcaption}
\usepackage{cleveref}
\usepackage[margin=1in]{geometry} 

\newtheorem{theorem}{Theorem}
\newtheorem{lemma}{Lemma}
\theoremstyle{definition} 

\newtheorem{remark}{Remark}
\newenvironment{subproof}[1][\proofname]{%
  \begin{proof}[#1]%
}{%
  \end{proof}%
}

\DeclareMathOperator{\intr}{int}
\DeclareMathOperator{\cl}{cl}
\DeclareMathOperator{\sym}{Sym}

\def\gif#1{\left\lfloor#1\right\rfloor} %Greatest Integer (Floor) Function e.g. \gif{x}
\newcommand{\vect}[2]{\left(\!\begin{smallmatrix} #1 \\ #2 \end{smallmatrix}\!\right)}
\def \N{\mathbbm{N}}
\def \Z{\mathbbm{Z}}

\def \R{\mathbbm{R}}

\makeatletter
\def\imod#1{\allowbreak\mkern10mu({\operator@font mod}\,\,#1)} %modulo symbol
\makeatother

\begin{document}

\captionsetup[subfigure]{subrefformat=simple,labelformat=simple,labelfont=rm}
\renewcommand\thesubfigure{(\alph{subfigure})}

\title[Fundamental Domains for Rhombic Lattices with Dihedral Symmetry of Order 8]{Fundamental Domains for Rhombic Lattices\\with Dihedral Symmetry of Order 8}
\author[J.R.C.G.~Damasco]{Joseph Ray Clarence G.~Damasco}
\address[J.R.C.G.~Damasco]{Institute of Mathematics, College of Science, University of the Philippines Diliman, Quezon City 1101, Philippines}
\email{jrcgdamasco@math.upd.edu.ph}

\author[D.~Frettl\"oh]{Dirk Frettl\"oh}
\address[D.~Frettl\"oh]{Technische Fakult\"at, Universit\"at Bielefeld, Postfach 100131, 33501 Germany}
\email{dfrettloeh@techfak.uni-bielefeld.de}

\author[M.J.C.~Loquias]{Manuel Joseph C.~Loquias}
\address[M.J.C.~Loquias]{Institute of Mathematics, College of Science, University of the Philippines Diliman, Quezon City 1101, Philippines}
\email{mjcloquias@math.upd.edu.ph}

\begin{abstract}
We show by construction that every rhombic lattice $\Gamma$ in $\mathbbm{R}^{2}$ 
has a fundamental domain whose symmetry group contains the point 
group of $\Gamma$ as a subgroup of index $2$. This solves the last open 
case of a question raised in \cite{df} on fundamental 
domains for planar lattices whose symmetry groups properly contain the 
point groups of the lattices.
\end{abstract}



\date{\today}

\maketitle

%%%%%%%%%%%%%%%%%%%%%%%%%%%%%%%%%%%%%%%%%%%%%%%%%%%%%%%%%%%%%%%%%%
\section{Introduction}

A {\em lattice} $\Gamma$ in $\R^d$ is the $\Z$-span of $d$ linearly independent
vectors in $\R^d$. The lattice $\Gamma$ forms a group that is isomorphic to the 
free abelian group of rank $d$. The point group $P(\Gamma)$ of $\Gamma$ is the set 
of Euclidean isometries in $\R^d$ fixing both $\Gamma$ and the origin. That is, 
$P(\Gamma)$ is a subgroup of the orthogonal group $\operatorname{O}(d)$  
consisting of those transformations that fix $\Gamma$. 
A {\em fundamental domain for} $\Gamma$ is a complete set of representatives of the 
orbits of $\R^{d}$ under the action of the group $\Gamma$. For example, 
$F=\big[\! -\tfrac{1}{2},\frac{1}{2}\big)^{2}$ is a fundamental domain for the lattice 
$\Gamma=\Z^{2} \subseteq \R^{2}$. We are interested in the geometric properties of 
fundamental domains, in particular, their symmetries. For a given $X \subseteq \R^{d}$,
the {\em symmetry group} $\sym(X)$ of $X$ is the group of Euclidean 
isometries in $\R^{d}$ that leave $X$ invariant. In our example, $\sym(F)$ consists of 
the identity and the reflection in the line $y=x$. Now, $F$ differs from the square region 
$Q=\big[\! -\tfrac{1}{2},\frac{1}{2}\big]^{2}$ only by a set of Lebesgue measure zero, 
and $Q$ is more symmetric compared to $F$ with $\sym(Q)=P(\Gamma)$. Note that throughout 
this paper, Lebesgue measure is used.

In order to discuss symmetries or geometric properties of fundamental domains, it is 
convenient to allow repetitions of representatives on sets of measure zero. In 
particular, we restrict ourselves to consider compact fundamental domains. Thus, in 
this paper, a fundamental domain for a lattice $\Gamma$ is a compact set $F$ such that 
(i) the Minkowski sum $F+\Gamma = \{p+\mathbf{v} \mid p \in F, \mathbf{v} \in \Gamma\}$ 
equals $\R^{d}$ and (ii) for any nonzero $\mathbf{v} \in \Gamma$, $\intr(F) \cap 
\intr(F+\mathbf{v}) = \varnothing$. 

The following fact generalizes the earlier example: For any lattice $\Gamma$, the 
{\em Voronoi region} $V$ of the origin is a fundamental domain for $\Gamma$ such 
that $\sym(V)$ equals $P(\Gamma)$ \cite{df}. Recall that the Voronoi region of a point 
$\mathbf{v}$ of a lattice $\Gamma$ in $\R^{d}$ is the collection of points in $\R^{d}$ 
closest to $\mathbf{v}$ than to any other point in $\Gamma$. The following question then 
naturally arises: Is there a fundamental domain $F$ for $\Gamma$ that has a symmetry group 
larger than $P(\Gamma)$, in the sense that $\sym(F)$ properly contains $P(\Gamma)$? 

About the year 2000 Veit Elser constructed a fundamental domain for $\Z^2$ 
exhibiting eightfold dihedral symmetry rather than the expected fourfold dihedral 
symmetry and a fundamental domain for the hexagonal lattice in $\R^2$ possessing 
twelvefold dihedral symmetry \cite{elserweb,elser}. 
Interestingly, this latter fundamental domain appears in different
contexts in \cite{bks} and \cite{coc}, where it serves as a window and an atomic 
surface, respectively, for mathematical quasicrystals. 

The procedures in \cite{elser} were refined and generalized in \cite{df} 
to obtain fundamental domains for rectangular lattices having the symmetries 
of a square. It was also shown that for oblique lattices, one can use a 
rectangular fundamental domain with one edge length equal to that of a 
basis vector, and a suitably chosen height. In fact, for a certain subclass of oblique 
lattices, this rectangular fundamental domain may be a square. Thus, if a bound for
$[\sym(F):P(\Gamma)]$ for planar lattices exists, it must be at least $4$.

Summarizing, for each planar lattice $\Gamma$
mentioned so far, there is a fundamental domain $F$ for $\Gamma$ for which 
$[\sym(F):P(\Gamma)]=2$ (Theorem 1.1 in \cite{df}). Hence, the only missing case for 
planar lattices is that of rhombic lattices and this paper is dedicated to 
fill this gap.

A {\em rhombic lattice} is a planar lattice for which there exists a basis 
consisting of two vectors of equal length such that the angle between the vectors 
is not $\frac{\pi}{2}$ (square lattice), or $\frac{\pi}{3}$ or $\frac{2\pi}{3}$ 
(hexagonal lattice). For any rhombic lattice $\Gamma$, $P(\Gamma)$ is isomorphic to 
the dihedral group $D_{2}$ of order $4$. In this paper, we construct for each rhombic 
lattice a fundamental domain $F$ with fourfold dihedral symmetry. Some of these 
fundamental domains are illustrated in Figures~\ref{fig:RhomIntRec} 
and~\ref{fig:RhomFDEx}. Our procedure differs substantially from the methods used in 
\cite{df}, as these methods do not seem to be applicable to rhombic lattices. Whereas 
the constructions in \cite{df} involve removing subregions and adding different ones 
at each step, the construction we employ here simply adds regions at each step. 

We prove our main result (Theorem~\ref{thm:rho}) as follows: First, we prove that it is 
enough to consider rhombic lattices whose diagonals have ratio at most $3$ 
(Lemma~\ref{lem:diagratio}). We then show that in order to construct $F$ it suffices 
to select a suitable subregion of a rectangle that satisfies certain conditions (see 
Rh1-Rh3 below). It turns out from Lemma~\ref{lem:expandrec} that it is enough to 
consider rectangles with edge ratio at most $3$. We then establish in 
Lemma~\ref{lem:fillpath} that after taking care of a particular subregion of the 
rectangle, the remaining portion of the rectangle may be dealt with by considering a 
smaller rectangle and finding a corresponding subregion satisfying analogues of 
conditions Rh1-Rh3, thereby creating an iterative process. Finally, we argue that the 
limit obtained by this iteration gives rise to a fundamental domain for $\Gamma$ 
satisfying the desired properties.

%%%%%%%%%%%%%%%%%%%%%%%%%%%%%%%%%%%%%%%%%%%%%%%%%%%%%%%%%%%%%%%%%%

\section{Main result}

\begin{theorem} \label{thm:rho}
Let $\Gamma$ be a rhombic lattice. Then there exists a fundamental domain $F$ for 
$\Gamma$ such that $[\sym(F):P(\Gamma)]=2$, in particular, $\sym(F) \cong D_{4}$.
\end{theorem}

\begin{proof}%[Proof of Theorem~\ref{thm:rho}] 
Let $\Gamma$ be a rhombic lattice. Because the properties we will be dealing with are 
invariant under similarities, we assume without loss of generality that $\Gamma$ has 
basis $\left\{\vect{2m}{0},\vect{m}{n}\right\}$, where $0<m<n$. Here, a basis for 
$\Gamma$ consisting of vectors of equal length is 
$\left\{\vect{m}{-n},\vect{m}{n}\right\}$. Let $\alpha$ be the reflection in the 
$y$-axis and $\beta$ be the half-turn about the origin (or central symmetry in $\R^{2}$). 
Then,  $P(\Gamma) = \langle \alpha, \beta \rangle \cong D_{2}$. 

Suppose first that $n \leq 3m$. Let $\mathcal{Q}$ be the triangular region with vertices 
$\vect{0}{0}$, $X=\vect{m}{0}$, and $W=\vect{0}{m}$. The union $\mathcal{S}$ of the 
images of $\mathcal{Q}$ under $P(\Gamma)$ forms a region with symmetry group $D_{4}$. 
Translates of $\mathcal{S}$ by distinct vectors in $\Gamma$ are either disjoint or 
intersect on sets of measure zero. See Figure~\ref{fig:RhomStepQ}.

\begin{figure}[ht]
\begin{center}
\includegraphics[width=.75\textwidth]{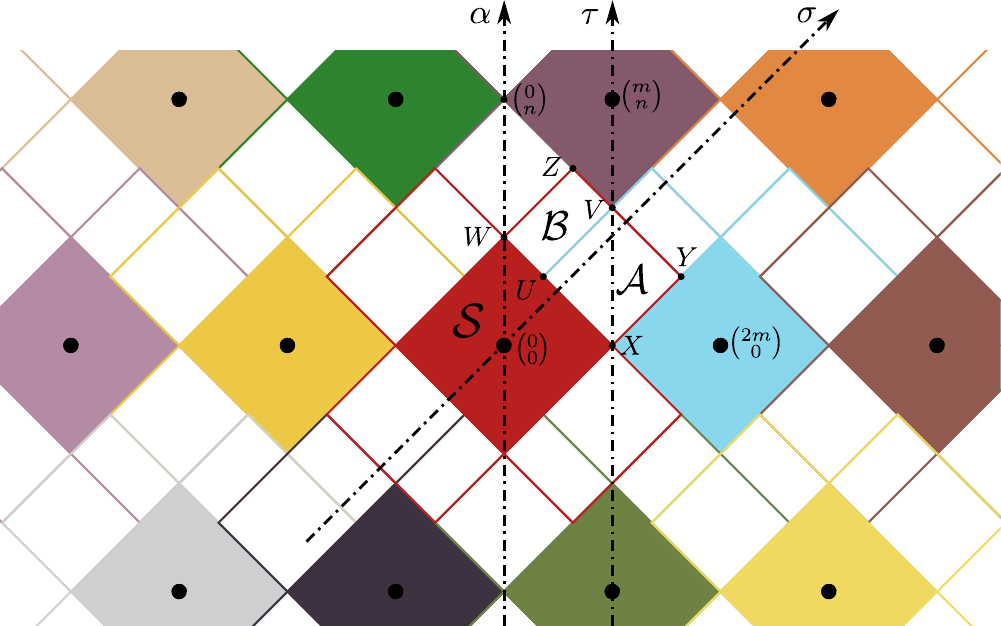}
\caption{First step in constructing a highly symmetric fundamental domain for a 
rhombic lattice.}
\label{fig:RhomStepQ}
\end{center}
\end{figure}

Now, consider the rectangular region $\mathcal{R}$ with vertices $W$, $X$, $Y$, and 
$Z$, where $YZ$ lies on the line through $\vect{0}{n}$ and $V=\vect{m}{n-m}$. We claim 
that there is a subset of $\mathcal{R}$ such that if $\mathcal{E}$ is the union of the 
images of this subset under $P(\Gamma)$, then the closure 
$F = \cl(\mathcal{S}\cup\mathcal{E})$ is a fundamental domain for $\Gamma$. If in 
addition the subset of $\mathcal{R}$ is preserved by the mirror reflection $\sigma$ in 
the perpendicular bisector of $WX$, then $F$ has symmetry group 
$\langle \alpha, \sigma \rangle \cong D_{4}$.

Since $n \leq 3m$, the ratio $\frac{WX}{XY} = \frac{2m}{(n-m)}$ is at least $1$. 
Partition $\mathcal{R}$ into the two sets: $\mathcal{A}$, the square $XYVU$, and 
$\mathcal{B} = \cl(\mathcal{R} \setminus \mathcal{A})$, the rectangle $WUVZ$ 
(see Figure~\ref{fig:RhomStepQ}). 
Note that if $n=3m$, then $V=Z$ and $\mathcal{B}$ is empty.

Consider the images of $\mathcal{A}$ and $\mathcal{B}$ under $P(\Gamma)$ and their 
translates by $\Gamma$. Equivalently, the union of these images and translates give 
$\sym(\Gamma)(\mathcal{A}\cup\mathcal{B})$. The only element of 
$\sym(\Gamma)\mathcal{A}$ that overlaps with $\mathcal{A}$ is 
$\alpha(\mathcal{A})+\vect{2m}{0}$. In fact, the two are equal. Similarly, 
$\mathcal{B} = \beta(\mathcal{B})+ \vect{m}{n}$, and $\mathcal{B}$ is disjoint with 
its other copies in $\sym(\Gamma)\mathcal{B}$. Thus, it suffices 
to find closed sets $\mathcal{K}\subseteq\mathcal{A}$ and 
$\mathcal{L}\subseteq\mathcal{B}$ such that 
$\mathcal{K} \cupdot \left(\alpha(\mathcal{K})+\vect{2m}{0}\right) = \mathcal{A}$, 
$\mathcal{L} \cupdot \left(\beta(\mathcal{L})+ \vect{m}{n}\right) = \mathcal{B}$, and 
$\mathcal{K} \cupdot \mathcal{L}$ is invariant under $\sigma$, where $\cupdot$ 
denotes a union of sets whose intersection has measure zero.  
We take the aforementioned $\mathcal{E}$ to be the union of the images of
 $\mathcal{K} \cupdot \mathcal{L}$ under $P(\Gamma)$.

The required conditions above can also be expressed in terms 
of symmetries of $\mathcal{A}$ and $\mathcal{B}$. Denote by $T_{\mathbf{v}}$ the 
translation by the vector $\mathbf{v}$. 
Let $\tau=T_{\vect{m}{0}}\alpha T^{-1}_{\vect{m}{0}}$ and 
$\rho=T_{\frac{1}{2}\vect{m}{n}}\beta T^{-1}_{\frac{1}{2}\vect{m}{n}}$. That is, 
$\tau$ is the mirror reflection in $XV$ and $\rho$ is the half-turn about 
the center of $\mathcal{B}$ which is $\frac{1}{2}\vect{m}{n}$. 
Using the fact that for any vector $\mathbf{v}$
and orthogonal linear transformation $\psi$, 
$\psi T_{\mathbf{v}} = T_{\psi(\mathbf{v})} \psi$ holds, we find that indeed,  
$\tau=T_{\vect{2m}{0}}\alpha$ and $\rho=T_{\vect{m}{n}}\beta$.
Thus, our main problem reduces to finding $\mathcal{K}\subseteq\mathcal{A}$ 
and $\mathcal{L}\subseteq\mathcal{B}$ satisfying the following properties.
	\begin{enumerate}[Rh1.]
	\item $\mathcal{K} \cupdot \tau(\mathcal{K}) = \mathcal{A}$
	\item $\mathcal{L} \cupdot \rho(\mathcal{L}) = \mathcal{B}$
	\item $\sigma(\mathcal{K}\cupdot\mathcal{L}) = \mathcal{K}\cupdot\mathcal{L}$
	\end{enumerate}
In the degenerate case $n=3m$, since $\mathcal{B}=\varnothing$, the problem reduces to 
finding $\mathcal{K}\subseteq\mathcal{A}$ that satisfies Rh1 and that is invariant 
under $\sigma$.

We now argue that the case where the ratio $\frac{n}{m}$ of the diagonals is greater 
than $3$ may be dealt with similarly.	

\begin{lemma}\label{lem:diagratio}
It suffices to consider the case when $n \leq 3m$.
\end{lemma}

\begin{subproof}[Proof of Lemma~\ref{lem:diagratio}]
If $n>3m$, let $t$ be the largest integer such that $n-2mt > m$. For each $i$ from 
$1$ to $t$, construct the square  $[m(i-1),mi] \times [m(i-1),mi]$. Obtain the union 
$\mathcal{C}$ of the images of these squares under $P(\Gamma)$. We note that 
$\sym(\mathcal{C}) \cong D_{4}$. By the choice of $t$, $\mathcal{C}$ does not intersect 
with any of its translates by non-horizontal vectors in $\Gamma$. See for example the 
squares in Figure~\ref{fig:RhomStepQ2a}, where $t=3$ and $\mathcal{C}$ is the union of 
the red squares. Note that the encircled dots represent lattice points. The points not 
encircled mark the furthest corners of the $\Gamma$-translates of $\mathcal{C}$.

\begin{figure}[ht]
\centering
\begin{subfigure}{0.47\textwidth}
\centering
\includegraphics[scale=0.7]{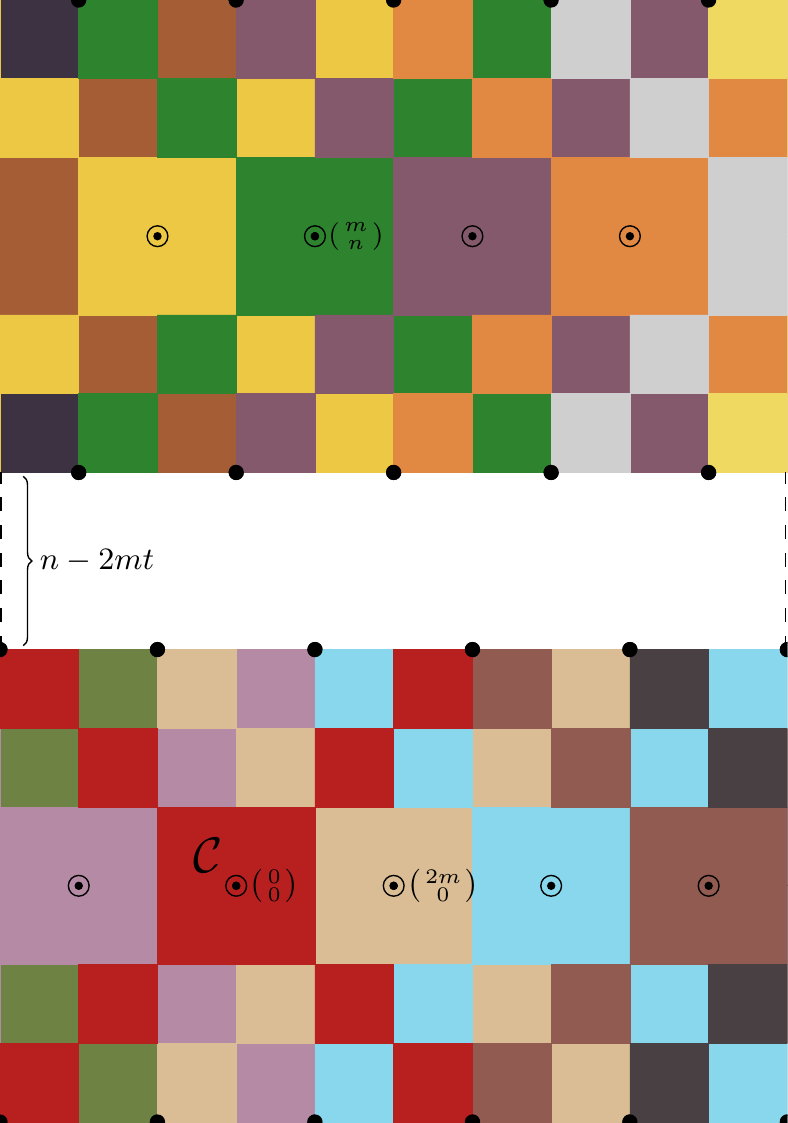}
\caption{}\label{fig:RhomStepQ2a}
\end{subfigure}
\begin{subfigure}{0.47\textwidth}
\centering
\includegraphics[scale=0.7]{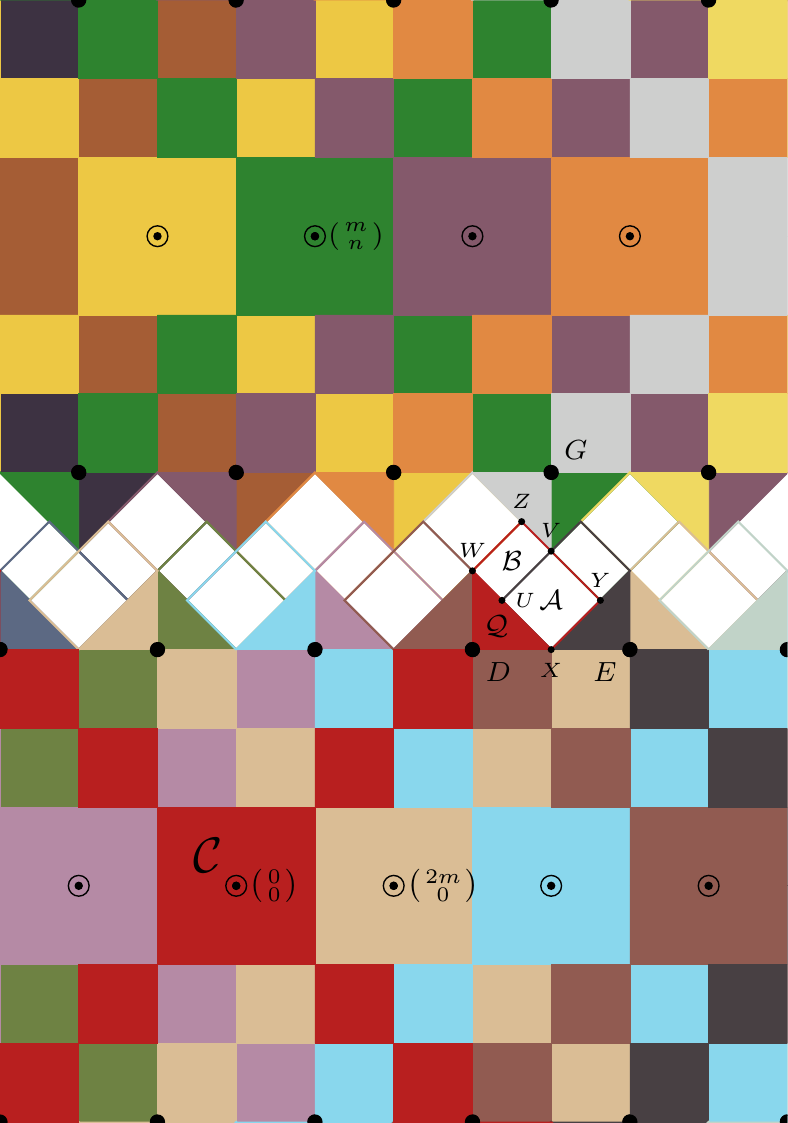}
\caption{}\label{fig:RhomStepQ2b}
\end{subfigure}
\caption{The case when the ratio of the diagonals is greater than $3$. Here, the 
smaller encircled dots represent lattice points while larger dots represent the 
furthest corners of copies of $\mathcal{C}$.}
\label{fig:RhomStepQ2}
\end{figure}

Furthermore, $\mathcal{C}$ does not overlap with its horizontal $\Gamma$-translates 
because the squares of the form $[m(i-1+2k),m(i+2k)] \times [m(i-1),mi]$, where 
$k \in \Z$, are precisely the horizontal $\Gamma$-translates of 
$[m(i-1),mi] \times [m(i-1),mi]$, while those of the form 
$[m(i+2k),m(i+1+2k)] \times [m(i-1),mi]$ are the horizontal $\Gamma$-translates of 
$[m(-i),m(-i+1)] \times [m(i-1),mi] = \alpha([m(i-1),mi] \times [m(i-1),mi])$. 
Thus, $\mathcal{C}$ does not overlap with its translates by nontrivial vectors in 
$\Gamma$. Moreover, the portion of the plane left uncovered by $\mathcal{C}+\Gamma$ 
is a union of horizontal strips of width $n-2mt$. 

Consider the points marked by dots on the boundary of one such horizontal strip, say 
the white strip in Figure~\ref{fig:RhomStepQ2a}. Each point on the lower boundary 
is both an upper-right corner and an upper-left corner of two horizontally-separated 
copies of $\mathcal{C}$, while each vertex on the upper boundary is both a lower-right 
and lower-left corner of two horizontally-separated copies of $\mathcal{C}$. 
The upper layer of points is the translation of the lower layer by $\vect{m}{n-2mt}$, 
and so the points on these two layers are spaced in the same way that the points of a 
rhombic lattice with basis $\vect{2m}{0}$ and $\vect{m}{n-2mt}$ are spaced. 
Note that by the choice of $t$, $m < n-2mt \leq 3m$.

In Figure~\ref{fig:RhomStepQ2b}, $D=\vect{mt}{mt}$ is the upper-right corner of 
$\mathcal{C}$, $E=\vect{m(t+2)}{mt}$ is the upper-left corner of 
$\mathcal{C} + (t+1)\vect{2m}{0}$, and $G=\vect{m(t+1)}{n-mt}$ is the lower-left 
corner of $\mathcal{C} + \vect{m}{n} + t\vect{2m}{0}$. This time, let $\mathcal{Q}$ be 
the triangular region with vertices $D=\vect{mt}{mt}$, $W = \vect{mt}{m(t+1)}$, and 
$X = \vect{m(t+1)}{mt}$. Since $m < n-2mt$, $\mathcal{Q}$ lies entirely in the strip. 
Let $\mathcal{S}$ be the union of the images of $\mathcal{Q}$ under $P(\Gamma)$. 
Extending the argument used for $\mathcal{C}$, $\mathcal{S}$ does not overlap with its 
horizontal translates. Furthermore, because $m < n-2mt$,  $\mathcal{S}$ does not 
overlap with its non-horizontal translates. Let $\mathcal{R}$ be the rectangle  
$WXYZ$ with $YZ$ on the line through $\vect{mt}{n-mt}$ and 
$V=\big(\begin{smallmatrix} m(t+1) \\ n-m(t+1) \end{smallmatrix}\big)$. Again, we can 
partition $\mathcal{R}$ into $\mathcal{A}$, the square $XYVU$, and 
$\mathcal{B} = \cl(\mathcal{R} \setminus \mathcal{A})$, the (possibly empty) rectangle 
$WUVZ$.  

Consider the images of $\mathcal{A}$ and $\mathcal{B}$ under $P(\Gamma)$ and the 
translates of these by $\Gamma$. Similar to the case when $n \leq 3m$, 
$\mathcal{A} = \alpha(\mathcal{A})+(t+1)\vect{2m}{0}$ and 
$\mathcal{B} = \beta(\mathcal{B})+ \vect{m}{n} + t\vect{2m}{0}$. Thus, if there exist 
$\mathcal{K}\subseteq \mathcal{A}$ and $\mathcal{L} \subseteq \mathcal{B}$ satisfying 
conditions Rh1, Rh2, and Rh3, and $\mathcal{E}$ is formed accordingly, then the set 
$F=\cl(\mathcal{C}\cup\mathcal{S}\cup\mathcal{E})$ is a fundamental domain for 
$\Gamma$ with symmetry group $D_{4}$.

This confirms that it suffices to consider the case when the ratio $\frac{n}{m}$ of the 
diagonals of the rhombic lattice is at most $3$. In fact, we use the same $\mathcal{Q}$, 
$\mathcal{K}$, and $\mathcal{L}$ (up to similarity) for two rhombic lattices if the 
difference between their diagonal ratios is an even integer.
\end{subproof}

We now construct the sets $\mathcal{K}$ and $\mathcal{L}$, where we assume 
$1 < \frac{n}{m} \leq 3$.  Let $a$ and $b$, with $a \geq b$, be the edge lengths of 
$\mathcal{R}$ so that $\frac{2m}{(n-m)} = \frac{a}{b}$. For convenience, 
we reorient $\mathcal{R}$ such that $U$ is at the origin, $\mathcal{A}=[0,b]\times[0,b]$, 
and $\mathcal{B}=[b-a,0]\times[0,b]$ if $a \neq b$. Recall that $\mathcal{B}$ is empty 
when $a=b$, which happens when $n=3m$. Recall also that $\sigma$ is the reflection in the 
perpendicular bisector of $WX$, $\tau$ is the reflection in the line through $XV$, and
$\rho$ is the half-turn about the center of $\mathcal{B}$. In the present orientation of 
$\mathcal{R}$ the axis of $\sigma$ is vertical and the axis of $\tau$ has slope $-1$.

\begin{lemma}\label{lem:intcase}
There exist $\mathcal{K}$ and $\mathcal{L}$ satisfying Rh1, Rh2, and Rh3 whenever 
the ratio $\frac{a}{b}$ is an integer.
\end{lemma}

\begin{subproof}[Proof of Lemma~\ref{lem:intcase}]
In each of the rectangles in 
Figure~\ref{fig:RhomIntRecStrips}, the rightmost square is $\mathcal{A}$. For 
$\frac{a}{b}=1$, $\mathcal{B}$ is empty and it suffices to find 
$\mathcal{K} \subseteq \mathcal{A} = \mathcal{R}$ such that 
$\mathcal{K} \cupdot \tau(\mathcal{K}) = \mathcal{A}$ and 
$\sigma(\mathcal{K})=\mathcal{K}$. One way to do this is to draw the two diagonals of 
$\mathcal{A}$ to divide it into four congruent triangles, and choose $\mathcal{K}$ to 
be the union of any two non-adjacent triangles, say the triangles marked red in 
Figure~\ref{fig:RhomIntRecStripsa}.

\begin{figure}[ht]
\centering
\begin{subfigure}{0.11\textwidth}
\centering
\includegraphics[scale=1]{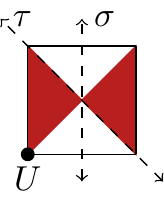}
\caption{}\label{fig:RhomIntRecStripsa}
\end{subfigure}
\begin{subfigure}{0.17\textwidth}
\centering
\includegraphics[scale=1]{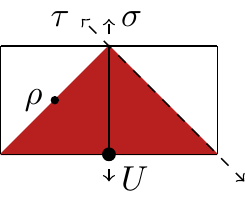}
\caption{}\label{fig:RhomIntRecStripsb}
\end{subfigure}
\begin{subfigure}{0.23\textwidth}
\centering
\includegraphics[scale=1]{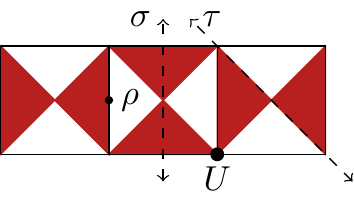}
\caption{}\label{fig:RhomIntRecStripsc}
\end{subfigure}
\begin{subfigure}{0.3\textwidth}
\centering
\includegraphics[scale=1]{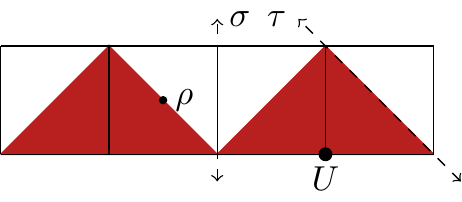}
\caption{}\label{fig:RhomIntRecStripsd}
\end{subfigure}
\caption{Choices for $\mathcal{K}$ and $\mathcal{L}$ when 
\subref{fig:RhomIntRecStripsa} $\frac{a}{b}=1$, 
\subref{fig:RhomIntRecStripsb} $\frac{a}{b}=2$, 
\subref{fig:RhomIntRecStripsc} $\frac{a}{b}=3$, 
\subref{fig:RhomIntRecStripsd} $\frac{a}{b}=4$.}
\label{fig:RhomIntRecStrips}
\end{figure}

If $\frac{a}{b}=2$, $\mathcal{B}$ is also a square, and $\mathcal{B}$ equals 
$\sigma({\mathcal{A}})$. This together with Rh3 implies that $\sigma(\mathcal{K})$ 
must be $\mathcal{L}$. Note that the chosen $\mathcal{K}$ for $\frac{a}{b}=1$ cannot 
be used, because its image under $\sigma$ violates Rh2. Rather, one can take 
$\mathcal{K}$ to be one of the triangles into which the axis 
of $\tau$ divides $\mathcal{A}$, as illustrated in Figure~\ref{fig:RhomIntRecStripsb}.

We now consider the case when $\frac{a}{b} \geq 3$. Because $\frac{a}{b}$ is an 
integer, we can divide $\mathcal{R}$ into a row of $\frac{a}{b}$ congruent squares each
of edge length $b$. If $\frac{a}{b}$ is odd, color triangles in each square as in the 
case when $\frac{a}{b}=1$, such that no two red triangles share a common edge 
(orientations of red triangles alternate). If $\frac{a}{b}$ is even, group the squares 
into adjacent pairs and color triangles in each pair as in the case when 
$\frac{a}{b}=2$. This is illustrated in Figures~\ref{fig:RhomIntRecStripsc} 
and~\ref{fig:RhomIntRecStripsd} for $\frac{a}{b}=3$ and $\frac{a}{b}=4$, respectively. 
In either case, if we let the red portion in $\mathcal{A}$ be $\mathcal{K}$ and the 
union of the other red portions be $\mathcal{L}$, then the sets $\mathcal{K}$ and 
$\mathcal{L}$ satisfy conditions Rh1, Rh2, and Rh3.
\end{subproof}

Figure~\ref{fig:RhomIntRec} shows the fundamental domains formed for the rhombic 
lattices with $\vect{m}{n} = \vect{3}{5}$ and $\vect{2}{3}$, which correspond to 
$\frac{a}{b}=3$ and $4$, respectively.

\begin{figure}[ht]
\centering
\begin{subfigure}{0.45\textwidth}
\centering
\includegraphics[scale=0.5]{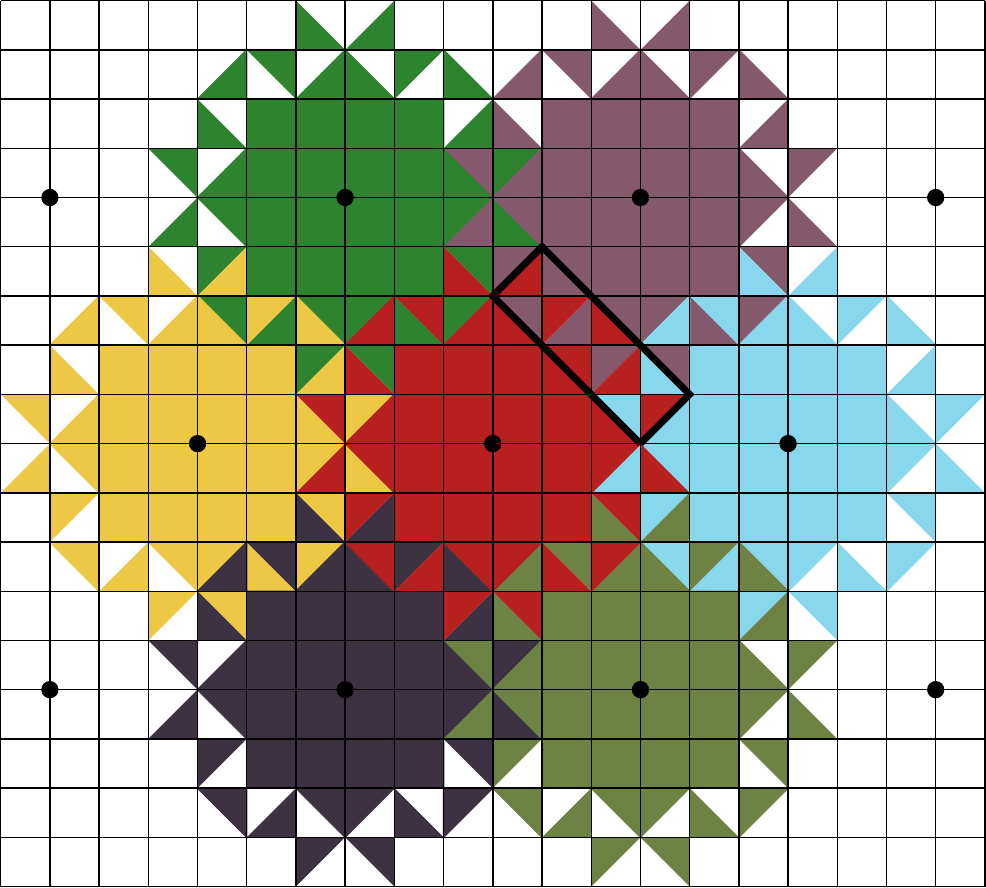}
\caption{}\label{fig:RhomIntReca}
\end{subfigure}
\begin{subfigure}{0.45\textwidth}
\centering
\includegraphics[scale=0.9]{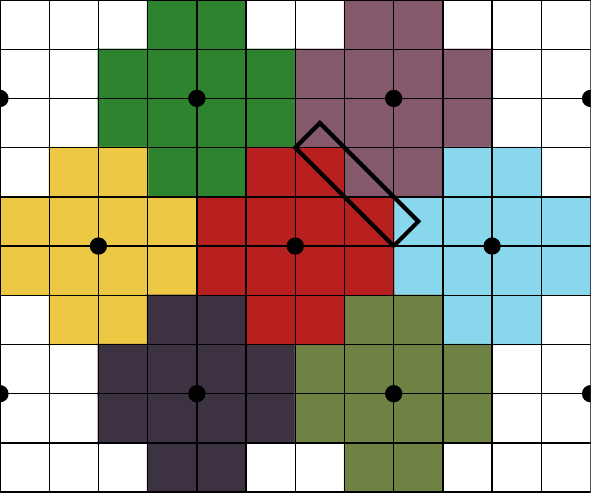}
\caption{}\label{fig:RhomIntRecb}
\end{subfigure}
\caption{Some fundamental domains where \subref{fig:RhomIntReca} $\frac{a}{b}=3$ and 
\subref{fig:RhomIntRecb} $\frac{a}{b}=4$.}
\label{fig:RhomIntRec}
\end{figure}

We shall use the same idea for more general rectangles $\mathcal{R}$, that is, we 
subdivide $\mathcal{R}$ into squares and take suitable subregions. If the rectangle 
has irrational edge ratio, the number of squares is necessarily infinite \cite{dehn}. 
The method we develop below works for any real value of the ratio, terminating after a 
finite number of steps precisely when the edge ratio is rational. First, we reduce the 
cases we have to consider.

\begin{lemma}
It suffices to construct $\mathcal{K}$ and $\mathcal{L}$ for rectangles $\mathcal{R}$ 
with edge ratio $\frac{a}{b}$ not exceeding $3$.
\label{lem:expandrec}
\end{lemma}

\begin{subproof}[Proof of Lemma~\ref{lem:expandrec}]
Suppose $\mathcal{R}=\mathcal{A} \cupdot \mathcal{B}$ has edge ratio 
$\frac{a}{b} \geq 1$, and $\mathcal{K}$ and $\mathcal{L}$ satisfy Rh1, Rh2, and Rh3. 
To prove the lemma, we show that for $\mathcal{R}' = \mathcal{A} \cupdot \mathcal{B}'$ 
with edge ratio $\frac{a}{b}+2$ and corresponding $\rho'$ and $\sigma'$, one can 
construct $\mathcal{L}' \subseteq \mathcal{B}'$ from $\mathcal{K}$ and $\mathcal{L}$ 
such that $\mathcal{K}$ and $\mathcal{L}'$ satisfy the corresponding analogues of Rh1, 
Rh2, and Rh3. 

Note that $\tau(\mathcal{A})=\mathcal{A}$, $\rho(\mathcal{B})=\mathcal{B}$, 
$\sigma(\mathcal{A}\cupdot \mathcal{B})=\mathcal{A}\cupdot \mathcal{B}$.
Refer to Figure~\ref{fig:RhomExtAB}. Consider the rectangle 
$\mathcal{R} \cupdot \rho(\mathcal{A}) =  \mathcal{A} \cupdot \mathcal{B} 
\cupdot \rho(\mathcal{A})$. We have 
\begin{equation*}
\rho(\mathcal{A} \cupdot \mathcal{B} \cupdot \rho(\mathcal{A}))  = \mathcal{A} 
\cupdot \mathcal{B} \cupdot \rho(\mathcal{A}).
\end{equation*}

\begin{figure}[ht]
\begin{center}
\includegraphics{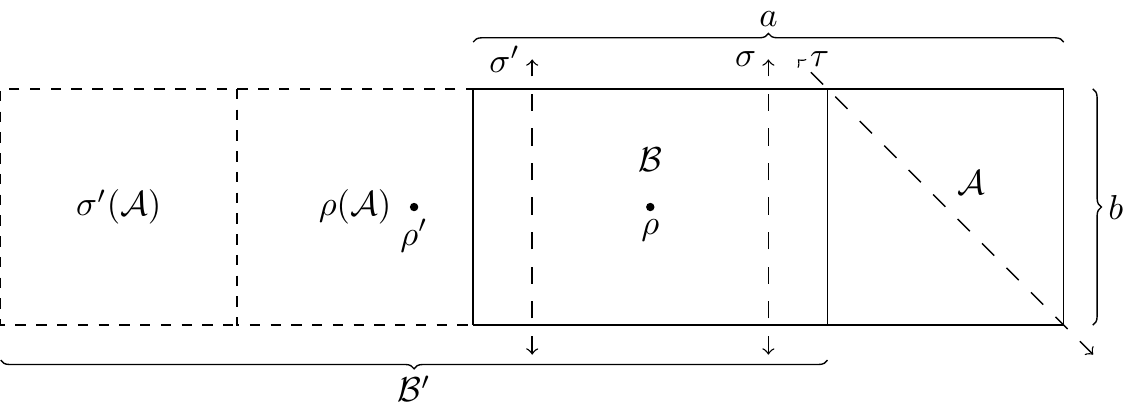} 
\caption{Constructing rectangle $\mathcal{R}'$ from rectangle $\mathcal{R}$ 
with edge ratio increased by $2$.}
\label{fig:RhomExtAB}
\end{center}
\end{figure}

Let $\sigma' = \rho \sigma \rho^{-1} = \rho \sigma \rho$, that is, $\sigma'$ is the 
reflection in the line formed by applying $\rho$ to the axis of $\sigma$. Then,
\begin{equation*}
\sigma' (\mathcal{B} \cupdot \rho (\mathcal{A})) = \rho \sigma \rho (\mathcal{B} 
\cupdot \rho (\mathcal{A})) = \rho \sigma (\mathcal{B} \cupdot \mathcal{A}) = 
\rho (\mathcal{B} \cupdot \mathcal{A}) =   \mathcal{B} \cupdot \rho (\mathcal{A}).
\end{equation*}

Set $\mathcal{B}' = \mathcal{B} \cup \rho(\mathcal{A}) \cup \sigma'(\mathcal{A})$  
and $\mathcal{R}'=\mathcal{A} \cupdot \mathcal{B}'$. Observe that 
the rectangle $\mathcal{R}'$ has edges with ratio $\frac{a}{b}+2$. Let $\rho' = 
\sigma' \rho (\sigma')^{-1} = \sigma' \rho \sigma'$, that is, $\rho'$ is
the half-turn about the image of the center of $\rho$ under $\sigma'$. It is easy
to verify that 
\begin{equation*}
\rho' (\mathcal{B}') = \mathcal{B}' \text{ and } \sigma' (\mathcal{R}') = \mathcal{R}'.
\end{equation*}

Set $\mathcal{L}' = \mathcal{L}\cupdot \rho\tau(\mathcal{K}) \cupdot 
\sigma'(\mathcal{K}) \subseteq \mathcal{B}'$. We claim that 
\begin{equation*} 
\sigma'(\mathcal{K}\cupdot\mathcal{L}')=\mathcal{K}\cupdot\mathcal{L}'
\text{ and }
\mathcal{L}'\cupdot\rho'(\mathcal{L}')=\mathcal{B}'.
\end{equation*} 

First, by Rh1 and Rh2, $\mathcal{L} \cupdot \rho \tau (\mathcal{K})$ and 
$\rho (\mathcal{L} \cupdot \mathcal{K})$ partition $\mathcal{B} \cupdot \rho 
(\mathcal{A})$, which is invariant under $\sigma'$. Because $\sigma' \rho 
(\mathcal{L} \cupdot \mathcal{K}) = \rho \sigma (\mathcal{L} \cupdot \mathcal{K}) 
= \rho (\mathcal{L} \cupdot \mathcal{K})$ by Rh3, we obtain that 
\begin{equation*}\sigma' (\mathcal{L} \cupdot \rho \tau (\mathcal{K})) = 
\mathcal{L} \cupdot \rho \tau (\mathcal{K}).
\end{equation*} 
It then follows that $\sigma' (\mathcal{K} \cupdot \mathcal{L}') = 
\mathcal{K} \cupdot \mathcal{L}'$.

Finally, we deduce that
\begin{align*}
\mathcal{L}' \cupdot \rho' (\mathcal{L}')
 & =  \mathcal{L}\cupdot \rho\tau(\mathcal{K}) \cupdot \sigma'(\mathcal{K}) \cupdot 
 \sigma' \rho \sigma'  (\mathcal{L}\cupdot \rho\tau(\mathcal{K}) \cupdot 
 \sigma'(\mathcal{K}))   \\
& = \sigma' (\mathcal{L}\cupdot \rho\tau(\mathcal{K})) \cupdot \sigma'(\mathcal{K}) 
\cupdot \sigma' \rho (\mathcal{L}\cupdot \rho\tau(\mathcal{K}) \cupdot \mathcal{K})   \\
& = \sigma' (\mathcal{A} \cupdot \mathcal{B} \cupdot \rho(\mathcal{A}))   
= \mathcal{B}'.
\end{align*}

Thus, for $\mathcal{R}'=\mathcal{A}\cupdot\mathcal{B}'$, $\tau$, $\rho'$, and 
$\sigma'$, and $\mathcal{K} \subseteq \mathcal{A}$ and 
$\mathcal{L}' \subseteq \mathcal{B}'$ as above, the corresponding analogues of 
conditions Rh1, Rh2, and Rh3 are satisfied. This shows how one can obtain from a 
given construction for a given rectangle a corresponding construction for a rectangle 
with edge ratio increased by 2. Hence, it does indeed suffice to consider those 
with edge ratio at most~3.
\end{subproof}

\begin{remark}
In Lemma~\ref{lem:expandrec}, the edge ratio $\frac{a}{b}$ may be an integer. The 
construction for the integer edge ratio case discussed in Lemma~\ref{lem:intcase} is 
precisely what is obtained when the procedure in the proof of Lemma~\ref{lem:expandrec} 
is applied inductively starting with the constructions for $\frac{a}{b} = 2$  and 
$\frac{a}{b}=3$.
\end{remark}

The reasoning in the proof of Lemma~\ref{lem:expandrec} applies as well to subsets of 
$\mathcal{A}$, $\mathcal{B}$, $\mathcal{K}$, and $\mathcal{L}$ that satisfy analogous 
conditions. We note them in the following remark as they will be useful later. 

\begin{remark} \ 
\begin{enumerate}[(a)]
\item Suppose $\mathcal{A}^{\ast} \subseteq \mathcal{A}$ and 
$\mathcal{B}^{\ast} \subseteq \mathcal{B}$ such that $\tau(\mathcal{A}^{\ast})
=\mathcal{A}^{\ast}$, $\rho(\mathcal{B}^{\ast})=\mathcal{B}^{\ast}$, 
$\sigma(\mathcal{A}^{\ast}\cupdot \mathcal{B}^{\ast})=\mathcal{A}^{\ast}\cupdot \mathcal{B}^{\ast}$. 
Then, $\mathcal{A}^{\ast}\cupdot \mathcal{B}^{\ast} \cupdot \rho(\mathcal{A}^{\ast}) \cupdot 
\sigma'(\mathcal{A}^{\ast})$ is invariant under $\sigma'$ and $\mathcal{B}^{\ast} \cupdot 
\rho(\mathcal{A}^{\ast}) \cupdot \sigma'(\mathcal{A}^{\ast})$ is invariant under $\rho'$.
\item Suppose $\mathcal{K}^{\ast} \subseteq \mathcal{A}^{\ast}$ and 
$\mathcal{L}^{\ast} \subseteq \mathcal{B}^{\ast}$ satisfy analogues of Rh1, Rh2, and Rh3, 
that is, $\mathcal{K}^{\ast} \cupdot \tau(\mathcal{K}^{\ast}) = \mathcal{A^{\ast}}$, 
$\mathcal{L}^{\ast} \cupdot \rho(\mathcal{L}^{\ast}) = \mathcal{B}^{\ast}$, $
\sigma(\mathcal{K^{\ast}}\cupdot\mathcal{L}^{\ast}) = \mathcal{K}^{\ast}\cupdot\mathcal{L}^{\ast}$. 
Then, $\mathcal{K}^{\ast} \cupdot \mathcal{L}^{\ast} \cupdot \rho\tau(\mathcal{K}^{\ast}) 
\cupdot  \sigma'(\mathcal{K}^{\ast})$ is invariant under $\sigma'$, and the region 
$\mathcal{L}^{\ast} \cupdot \rho\tau(\mathcal{K}^{\ast}) \cupdot \sigma'(\mathcal{K}^{\ast})$  
and its image under $\rho'$ partition 
$\mathcal{B}^{\ast} \cupdot \rho(\mathcal{A}^{\ast}) \cupdot \sigma'(\mathcal{A}^{\ast})$.
\end{enumerate}
\label{rem:expandrec}
\end{remark}

We have seen how one can directly construct $\mathcal{K}$ and $\mathcal{L}$ when 
$\frac{a}{b}$ is an integer. For non-integer values of $\frac{a}{b}$, we will construct 
$\mathcal{K}$ and $\mathcal{L}$ in a finite or countable number of steps. In particular, we will 
show that there are closed regions $\mathcal{A}_{i}^{\ast}$, $\mathcal{B}_{i}^{\ast}$, 
$\mathcal{K}_{i}^{\ast}$, $\mathcal{L}_{i}^{\ast}$ for $i \in \mathbbm{N}$ 
with the following properties. 

\begin{enumerate}[D1.]
\item For every $i$,
	\begin{enumerate}[a.]
	\item $\tau(\mathcal{A}_{i}^{\ast})=\mathcal{A}_{i}^{\ast}$, 
	$\rho(\mathcal{B}_{i}^{\ast})=\mathcal{B}_{i}^{\ast}$, $\sigma(\mathcal{A}_{i}^{\ast}\cupdot 
	\mathcal{B}_{i}^{\ast})=\mathcal{A}_{i}^{\ast}\cupdot\mathcal{B}_{i}^{\ast}$.
	\item $\mathcal{K}_{i}^{\ast}\subseteq \mathcal{A}_{i}^{\ast}$ and $\mathcal{L}_{i} ^{\ast}
	\subseteq 	\mathcal{B}_{i}^{\ast}$, and they satisfy analogues of Rh1, Rh2, and Rh3. 
	\end{enumerate}
\item $\mathcal{A}=\cl\left(\displaystyle\bigcupdot_{i}\mathcal{A}_{i}^{\ast}\right)$ and
	$\mathcal{B}=\cl\left(\displaystyle\bigcupdot_{i}\mathcal{B}_{i}^{\ast}\right)$. 
\end{enumerate}

If so, let $\mathcal{K} \coloneqq \cl\left(\displaystyle\bigcupdot_{i}\mathcal{K}_{i}^{\ast}\right)
\subseteq \mathcal{A}$, and $\mathcal{L} \coloneqq 
\cl\left(\displaystyle\bigcupdot_{i}\mathcal{L}_{i}^{\ast}\right) \subseteq \mathcal{B}$. Then if
in addition, the unions $\mathcal{K} \cup \tau(\mathcal{K})$ and $\mathcal{L} \cup \rho(\mathcal{L})$
are disjoint up to sets of measure zero, then $\mathcal{K}$ and $\mathcal{L}$ satisfy
Rh1, Rh2, and Rh3.

Set $a_{1}=a$, $b_{1}=b$. Suppose $1<\frac{a_{1}}{b_{1}}<2$, that is, $\mathcal{A}$ is
larger than $\mathcal{B}$. Let $v_{1}= b_{1}-(a_{1}-b_{1})
= 2b_{1}-a_{1}$ and $q_{1}=\gif{\frac{b_{1}}{v_{1}}}$. Recall that 
$\mathcal{A}=[0,b]\times[0,b]$, and $\mathcal{B}=[b-a,0]\times[0,b]$. Consider all the 
squares of the form $[v_{1}i_{1},v_{1}(i_{1}+1)] \times 
[b_{1}-v_{1}(j_{1}+1),b_{1}-v_{1}j_{1}]$, where $i_{1}$ and $j_{1}$ are nonnegative 
integers such that $i_{1}+j_{1} \leq  q_{1}-1$. See for example 
Figure~\ref{fig:RhomA1B112}, where $q_{1}=3$.  The squares all have edge length $v_{1}$, 
and they form a triangular array of squares propped against the left and upper edges of 
$\mathcal{A}$. We say that the array is anchored at the upper-left vertex of 
$\mathcal{A}$, and is of order $q_{1}$, with edge length $v_{1}$. If in addition, 
$q_{1}>1$, construct a second array of squares anchored at the lower-right vertex of 
$\mathcal{A}$ of order $q_{1}-1$ with edge length $v_{1}$. This is the largest order of 
a second array that can be constructed such that the two arrays do not overlap. The union 
$\mathcal{A}_{1}^{\ast}$ of the two arrays is invariant under $\tau$. Next, form two arrays 
within $\mathcal{B}$, both of order $q_{1}-1$ and edge length $v_{1}$, one anchored at 
the upper-right vertex, and another at the lower-left vertex of 
$\mathcal{B}$. Let $\mathcal{B}_{1}^{\ast}$ be the union of these two arrays. Then 
$\rho(\mathcal{B}_{1}^{\ast})=\mathcal{B}_{1}^{\ast}$ and 
$\sigma(\mathcal{A}_{1}^{\ast}\cupdot \mathcal{B}_{1}^{\ast}) = 
\mathcal{A}_{1}^{\ast}\cupdot\mathcal{B}_{1}^{\ast}$, so that $\mathcal{A}_{1}^{\ast}$ and 
$\mathcal{B}_{1}^{\ast}$ satisfy D1a.

\begin{figure}[ht]
\centering
\begin{subfigure}{0.34\textwidth}
\centering
\includegraphics[scale=0.58]{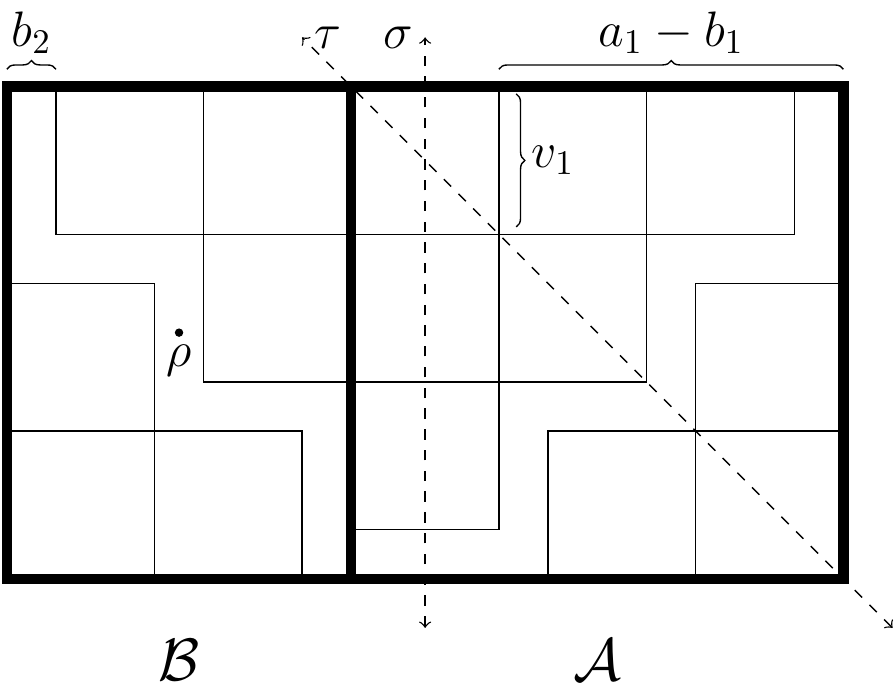}
\caption{}\label{fig:RhomA1B112}
\end{subfigure}
\begin{subfigure}{0.48\textwidth}
\centering
\includegraphics[scale=0.58]{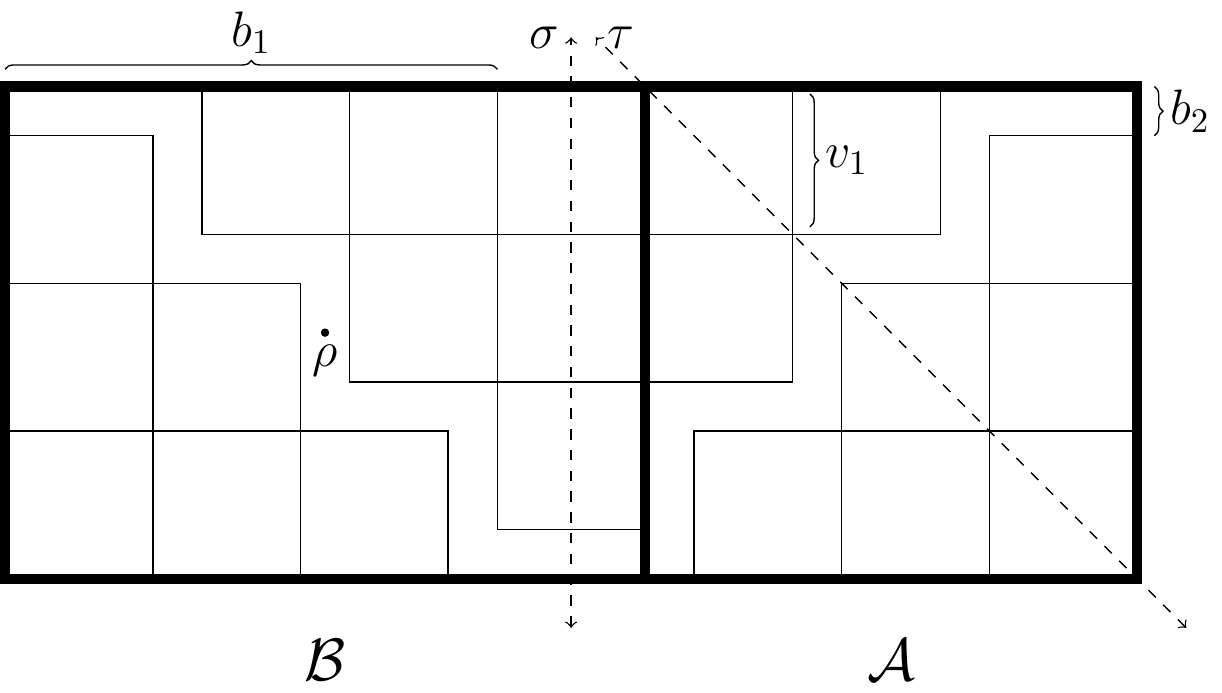}
\caption{}\label{fig:RhomA1B123}
\end{subfigure}
\caption{First steps in constructing $\mathcal{K}$ and $\mathcal{L}$ for 
\subref{fig:RhomA1B112} $1<\frac{a}{b}<2$ and \subref{fig:RhomA1B123} 
$2<\frac{a}{b}<3$.}
\label{fig:RhomA1B1}
\end{figure}

Suppose $2<\frac{a_{1}}{b_{1}}<3$. This time, $\mathcal{B}$ is larger than $\mathcal{A}$ 
and we set $v_{1}=(a_{1}-b_{1})-b_{1}=a_{1}-2b_{1}$ and $q_{1}=\gif{\frac{b_{1}}{v_{1}}}$. 
See Figure~\ref{fig:RhomA1B123}. Within $\mathcal{A}$, construct two arrays of 
edge length $v_{1}$, one with order $q_{1}$  anchored at the lower-right vertex of 
$\mathcal{A}$, and the other with order $q_{1}-1$ (if positive) anchored at the upper-left 
vertex of $\mathcal{A}$. We denote the union of these arrays by $\mathcal{A}_{1}^{\ast}$. 
Then form two arrays within $\mathcal{B}$, both of order $q_{1}$ and edge length $v_{1}$, 
one anchored at the upper-right vertex of $\mathcal{B}$, and another at the lower-left 
vertex of $\mathcal{B}$. Let $\mathcal{B}_{1}^{\ast}$ be the union of these two arrays.  
Again, D1a is satisfied.

To construct $\mathcal{K}_{1}^{\ast}$ and $\mathcal{L}_{1}^{\ast}$, draw the diagonals of each 
square in each array, thereby dividing each square into four congruent triangles  
(as done in the odd-integer-ratio case). Consider the arrays in $\mathcal{A}_{1}$ and 
$\mathcal{B}_{1}$ anchored at a shared vertex of $\mathcal{A}$ and $\mathcal{B}$. 
For each square in these arrays, color the upper and lower triangles red. For the 
remaining two arrays, color the left and right triangles 
of each square red as in Figure~\ref{fig:RhomK1L1}. Let $\mathcal{K}_{1}^{\ast}$  and 
$\mathcal{L}_{1}^{\ast}$ be the union of the red regions contained in $\mathcal{A}_{1}^{\ast}$ 
and $\mathcal{B}_{1}^{\ast}$, respectively. By construction, 
$\mathcal{K}_{1}^{\ast}$ and $\mathcal{L}_{1}^{\ast}$ satisfy D1b.

\begin{figure}[ht]
\centering
\begin{subfigure}{0.34\textwidth}
\centering
\includegraphics[scale=0.58]{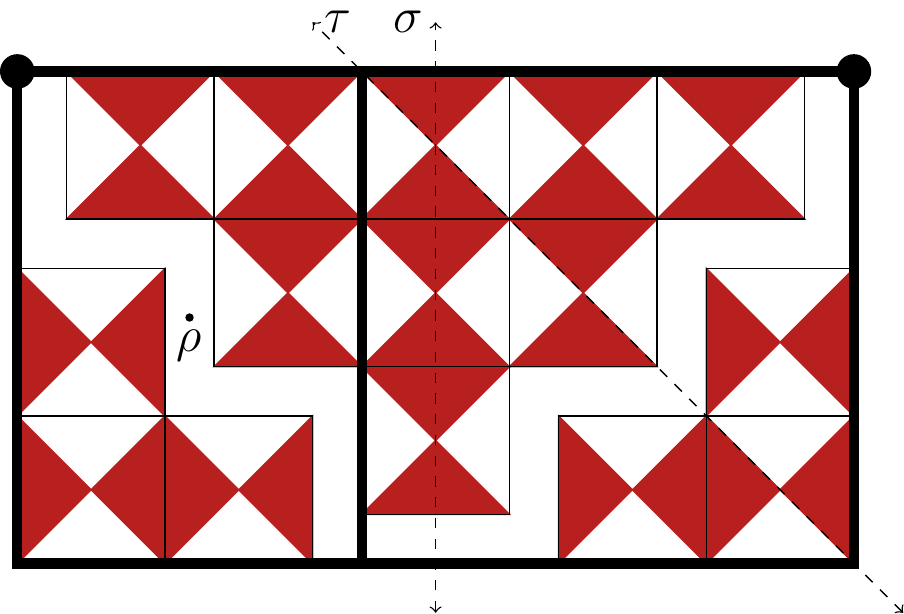}
\caption{}\label{fig:RhomK1L112}
\end{subfigure}
\begin{subfigure}{0.45\textwidth}
\centering
\includegraphics[scale=0.58]{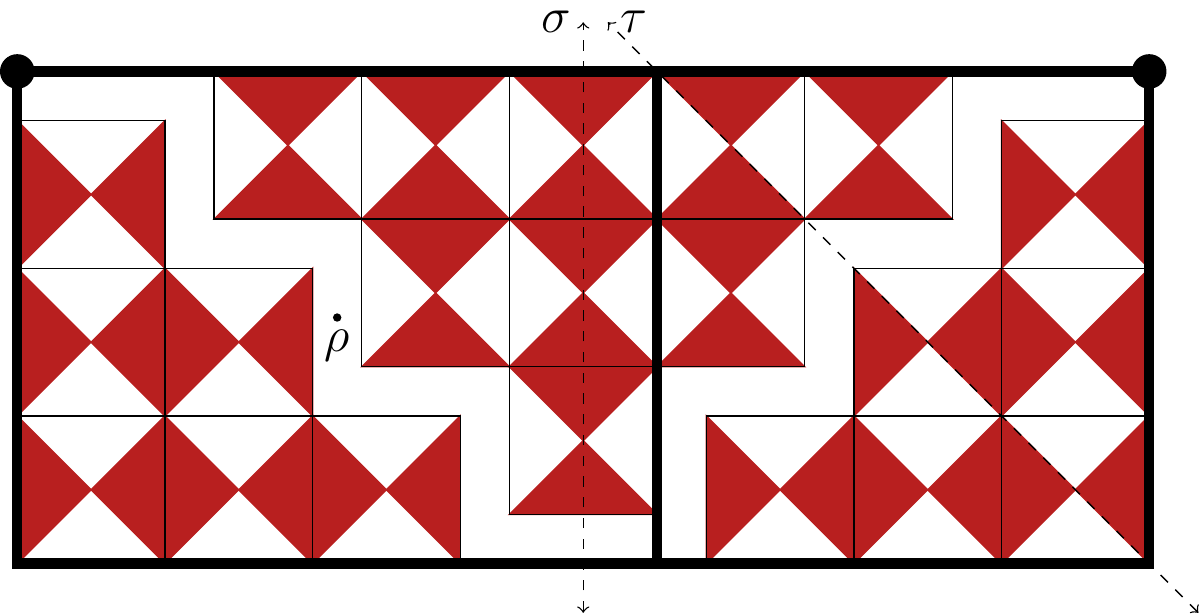}
\caption{}\label{fig:RhomK1L123}
\end{subfigure}
\caption{Subregions $\mathcal{K}_{1}^{\ast}$ and $\mathcal{L}_{1}^{\ast}$ of 
$\mathcal{A}_{1}^{\ast}$ and $\mathcal{B}_{1}^{\ast}$ for \subref{fig:RhomK1L112} 
$1<\frac{a}{b}<2$ and \subref{fig:RhomK1L123} $2<\frac{a}{b}<3$.}
\label{fig:RhomK1L1}
\end{figure}

Let $b_{2}=b_{1}-v_{1}q_{1}$. This is the uniform width of the ``path'' $\mathcal{S}_{1}$ 
left uncovered by $\mathcal{A}_{1}^{\ast}$ and $\mathcal{B}_{1}^{\ast}$ in $\mathcal{R}$. 
Thus, if $\frac{b_{1}}{v_{1}}$ is an integer, then $\mathcal{A}=\mathcal{A}_{1}^{\ast}$ 
and $\mathcal{B}=\mathcal{B}_{1}^{\ast}$, and so we are done if we set 
$\mathcal{K}=\mathcal{K}_{1}^{\ast}$ and $\mathcal{L}=\mathcal{L}_{1}^{\ast}$. Otherwise, 
consider the central column of squares, that is, the column fixed by $\sigma$. The squares 
in this column cover a vertical strip of $\mathcal{R}$ of width $v_{1}$ except for a 
rectangle $\hat{\mathcal{R}}$ with dimensions $b_{2} \times v_{1}$, where $b_{2}<v_{1}$. 
Partition $\hat{\mathcal{R}}$ into a square $\hat{\mathcal{A}}$ of edge length $b_{2}$ 
and a rectangle $\hat{\mathcal{B}}$ as in Figure~\ref{fig:RhomA2B2Zoom}. Let $\hat{\rho}$ 
be the half-turn about the center of $\hat{\mathcal{B}}$ and $\hat{\tau}$ be the reflection 
in the line containing the diagonal of $\hat{\mathcal{A}}$ that passes through one corner 
of the central column of squares. 

\begin{figure}[ht]
\begin{center}
\includegraphics[scale=0.85]{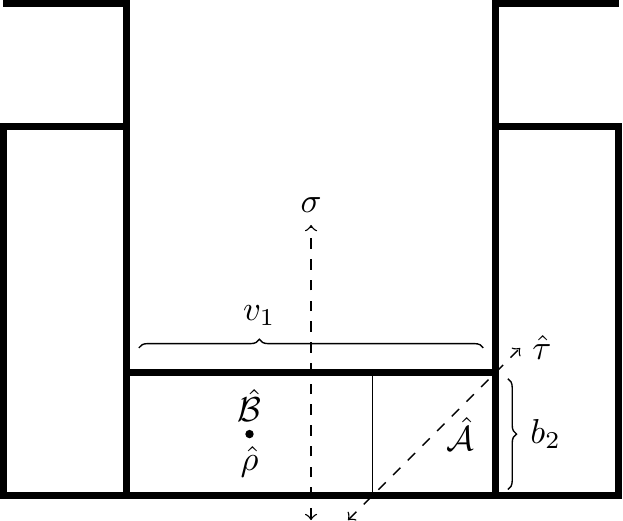}
\caption{The rectangle $\hat{\mathcal{R}} = \hat{\mathcal{A}} \cupdot \hat{\mathcal{B}}$.}
\label{fig:RhomA2B2Zoom}
\end{center}
\end{figure}

\begin{lemma}
It suffices to find $\hat{\mathcal{K}}\subseteq \hat{\mathcal{A}}$ and $\hat{\mathcal{L}} 
\subseteq \hat{\mathcal{B}}$ that satisfy the corresponding analogues of Rh1, Rh2, and Rh3.
\label{lem:fillpath}
\end{lemma}

\begin{subproof}[Proof of Lemma~\ref{lem:fillpath}]

Suppose $1<\frac{a_{1}}{b_{1}}<2$. Then, $\hat{\mathcal{R}} \subseteq \mathcal{A}$. 
Consider the sequences $\left\{R_{j}\right\}$ and $\left\{F_{j}\right\}$ defined by 
$R_{1} = \hat{\mathcal{R}} = \hat{\mathcal{A}} \cupdot \hat{\mathcal{B}}$, 
$F_{1}=\hat{\mathcal{K}} \cupdot \hat{\mathcal{L}}$ , and 
for $1 \leq j \leq 4q_{1}$, 
$$R_{j+1}=\left\{\begin{aligned} 
{\tau}(R_{j}), & && j \equiv 1 \imod 4\\
\sigma(R_{j}), & && j \equiv 0 \imod 2\\
{\rho}(R_{j}), & && j \equiv 3 \imod 4
\end{aligned}\right. \ \text{ and } \ 
F_{j+1}=\left\{\begin{aligned} 
{\tau}(\cl(R_{j}\setminus F_{j})), & && j \equiv 1 \imod 4\\
\sigma(F_{j}), & && j \equiv 0 \imod 2\\
{\rho}(\cl(R_{j}\setminus F_{j})), & && j \equiv 3 \imod 4.
\end{aligned}\right.
$$
The rectangles ${R}_{j}$ generated are labeled in Figure~\ref{fig:RhomA2B2} by their 
indices. Note that the orientations of the last rectangles in the sequence are not 
necessarily the same as those in the figure. The longer edges of $R_{4q_{1}-3}$ and 
$R_{4q_{1}-2}$ may be horizontal and vertical, respectively. Regardless of the 
orientations,the following properties hold.

\begin{enumerate}[(a)]
\item For $0 \leq k \leq q_{1}-1$, $R_{4k+1}$ and $R_{4k+2}$ are non-overlapping 
$\tau$-images of each other, and $F_{4k+1} \cupdot F_{4k+2}$ and its image under 
$\tau$ partition $R_{4k+1} \cupdot R_{4k+2}$.
\item For $0 \leq k \leq q_{1}-2$, $R_{4k+3}$ and $R_{4k+4}$ are non-overlapping 
$\rho$-images of each other, and $F_{4k+3} \cupdot F_{4k+4}$ and its image under 
$\rho$ partition $R_{4k+3} \cupdot R_{4k+4}$.
\item The regions $R_{4q_{1}-1}$ and $R_{4q_{1}}$ overlap on a rectangle congruent 
to $\hat{\mathcal{B}}$. By the choice of $\hat{\mathcal{L}}$, the portions of 
$F_{4q_{1}-1}$ and $F_{4q_{1}}$ inside this rectangle coincide completely.
\item For $1 \leq k \leq 2q_{1}$, $R_{2k}$ and $R_{2k+1}$ are $\sigma$-images of 
each other. The same holds true for $F_{2k}$ and $F_{2k+1}$.
\item The union $\displaystyle\bigcup_{j=1}^{4q_{1}} R_{j}$ covers 
$\mathcal{R}\setminus(\mathcal{A}_{1}^{\ast}\cup\mathcal{B}_{1}^{\ast})$ except for a 
square region $\mathcal{H}$. The rectangle $R_{4q_{1}+1}$ covers this. By the choice of 
$\hat{\mathcal{K}}$, the portion of $F_{4q_{1}+1}$ inside this square and the 
image of this portion under $\tau$ partition the square. 
\end{enumerate}

\begin{figure}[ht]
\begin{center}
\includegraphics[scale=0.8]{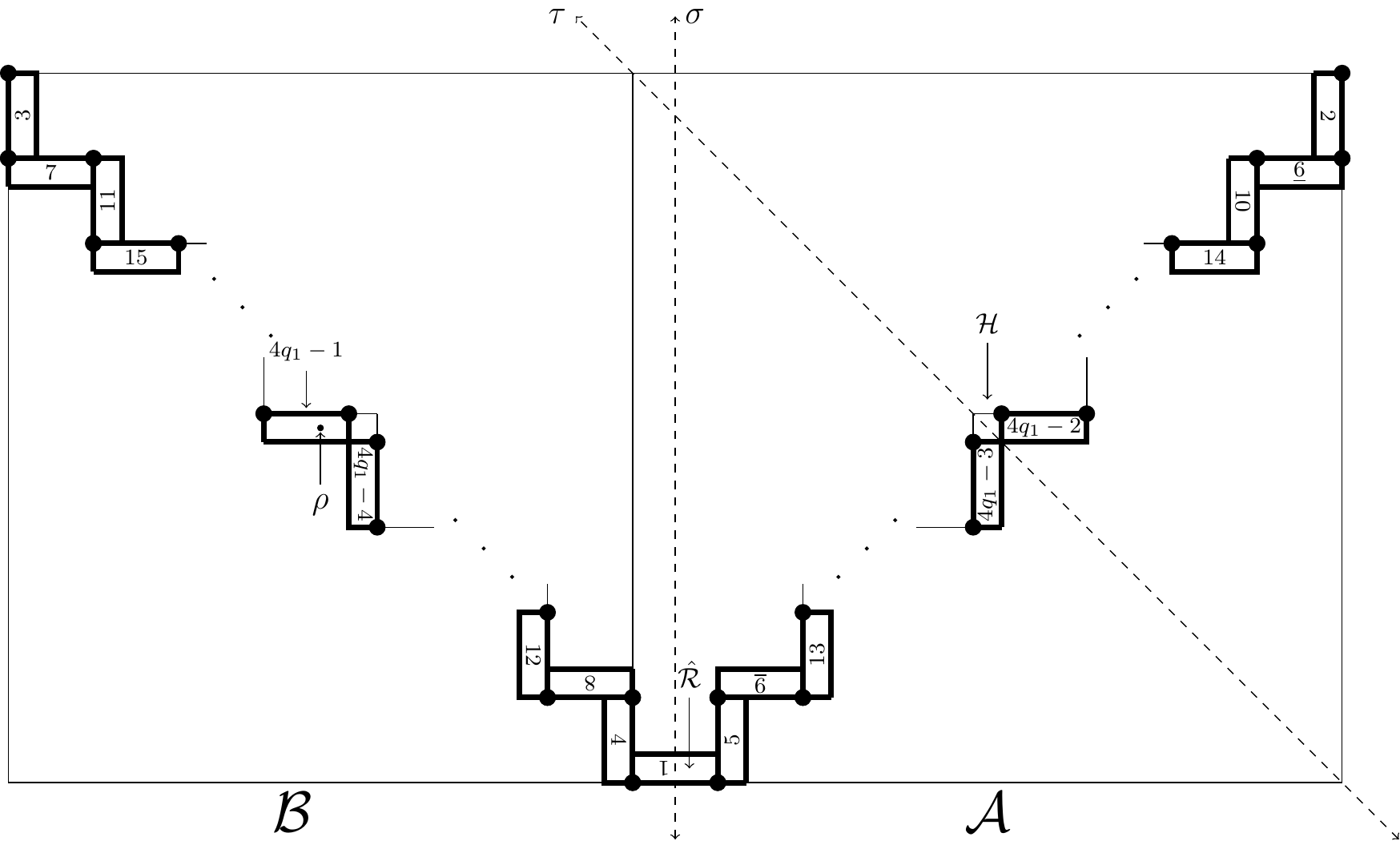}
\caption{Filling up the uncovered path $\mathcal{S}_{1}$ in $\mathcal{R}$ by applying
a sequence of $\sigma$, $\tau$, and $\rho$ on $\hat{\mathcal{R}}$. The two vertices of 
rectangle $R_{i}$ marked by big dots are the corresponding images of the two vertices 
marked in $R_{i-1}$.}
\label{fig:RhomA2B2}
\end{center}
\end{figure}

Define $\mathcal{A}_{2}^{\ast} = \displaystyle \bigcup_{j \equiv 1,2 \imod 4} R_{j} 
\subseteq \mathcal{A}$,  $\mathcal{K}_{2}^{\ast} = \displaystyle \bigcup_{j \equiv 1,2 \imod 4} 
F_{j} \subseteq \mathcal{A}_{2}^{\ast}$, $\mathcal{B}_{2}^{\ast} = \displaystyle 
\bigcup_{j \equiv 0,3 \imod 4} R_{j} \subseteq \mathcal{B}$, and $\mathcal{L}_{2}^{\ast} =
\displaystyle \bigcup_{j \equiv 0,3 \imod 4} F_{j} \subseteq \mathcal{B}_{2}^{\ast}$. We 
see that these, with $\tau$, $\rho$, and $\sigma$, satisfy analogues of Rh1, Rh2, 
and Rh3. Moreover, $\mathcal{A}_{i}^{\ast}$, $\mathcal{B}_{i}^{\ast}$, $\mathcal{K}_{i}^{\ast}$, 
and $\mathcal{L}_{i}^{\ast}$ for $i = 1, 2$ satisfy conditions D1 and D2.

The case when $2<\frac{a_{1}}{b_{1}}<3$ is treated analogously. In this case, 
$\hat{\mathcal{R}} \subseteq \mathcal{B}$. Define the sequences $\left\{R_{j}\right\}$ 
and $\left\{F_{j}\right\}$ by $R_{1} = \hat{\mathcal{R}} = 
\hat{\mathcal{A}} \cupdot \hat{\mathcal{B}}$, $F_{1}=\hat{\mathcal{K}} \cupdot 
\hat{\mathcal{L}}$, and for $1 \leq j \leq 4q_{1}+2$, 
$$R_{j+1}=\left\{\begin{aligned} 
{\rho}(R_{j}), & && j \equiv 1 \imod 4\\
\sigma(R_{j}), & && j \equiv 0 \imod 2\\
{\tau}(R_{j}), & && j \equiv 3 \imod 4
\end{aligned}\right. \ \text{ and } \ 
F_{j+1}=\left\{\begin{aligned} 
{\rho}(\cl(R_{j}\setminus F_{j})), & && j \equiv 1 \imod 4\\
\sigma(F_{j}), & && j \equiv 0 \imod 2\\
{\tau}(\cl(R_{j}\setminus F_{j})), & && j \equiv 3 \imod 4.
\end{aligned}\right.
$$

This time, take $\mathcal{A}_{2}^{\ast} = \displaystyle \bigcup_{j \equiv 0,3 \imod 4} R_{j} 
\subseteq \mathcal{A}$,  $\mathcal{K}_{2}^{\ast} = \displaystyle \bigcup_{j \equiv 0,3 \imod 4} 
F_{j} \subseteq \mathcal{A}_{2}^{\ast}$, $\mathcal{B}_{2}^{\ast} = \displaystyle 
\bigcup_{j \equiv 1,2 \imod 4} R_{j} \subseteq \mathcal{B}$, and $\mathcal{L}_{2}^{\ast} =
\displaystyle \bigcup_{j \equiv 1,2 \imod 4} F_{j} \subseteq \mathcal{B}_{2}^{\ast}$.

\end{subproof}

\begin{remark} \ 
As in Lemma~\ref{lem:expandrec}, the properties enumerated in the proof of Lemma 
\ref{lem:fillpath} are true for suitably chosen subsets of $\hat{\mathcal{A}}$
and $\hat{\mathcal{B}}$. In particular, suppose
$\hat{\mathcal{A}}^{\ast} \subseteq \hat{\mathcal{A}}$ and $\hat{\mathcal{B}}^{\ast} \subseteq 
\hat{\mathcal{B}}$ such that $\hat{\tau}(\hat{\mathcal{A}}^{\ast})=\hat{\mathcal{A}}^{\ast}$, 
$\hat{\rho}(\hat{\mathcal{B}}^{\ast})=\hat{\mathcal{B}}^{\ast}$, $\sigma(\hat{\mathcal{A}}^{\ast}\cupdot 
\hat{\mathcal{B}}^{\ast})=\hat{\mathcal{A}}^{\ast}\cupdot \hat{\mathcal{B}}^{\ast}$. Suppose 
$\hat{\mathcal{K}}^{\ast} \subseteq \hat{\mathcal{A}}^{\ast}$ and $\hat{\mathcal{L}}^{\ast} 
\subseteq \hat{\mathcal{B}}^{\ast}$ satisfy analogues of Rh1, Rh2, and Rh3. Obtain the sequences 
$\{{R}_{j}^{\ast}\}$ and  $\{{F}_{j}^{\ast}\}$ analogously, and form $\hat{\mathcal{A}}_{2}^{\ast}$, 
$\hat{\mathcal{B}}_{2}^{\ast}$, $\hat{\mathcal{K}}_{2}^{\ast}$, $\hat{\mathcal{L}}_{2}^{\ast}$
accordingly. Then, these also satisfy analogues of Rh1, Rh2, and Rh3.
\label{rem:fillpath}
\end{remark}
	
We thus have the following procedure.\smallskip

\noindent\textbf{Step 1:} \ 
	\begin{itemize}
	\item Without loss of generality, assume that the rectangle $\mathcal{R}_{1} \coloneqq 
	\mathcal{R}$ has a non-integer edge ratio between $1$ and $3$. Otherwise we can use Lemma 
	\ref{lem:intcase} or Lemma~\ref{lem:expandrec}.
	\item Fill $\mathcal{R}_{1}$ with $\mathcal{A}_{1}^{\ast}$, $\mathcal{B}_{1}^{\ast}$, 
	$\mathcal{K}_{1}^{\ast}$, and $\mathcal{L}_{1}^{\ast}$ as described above, leaving the 
	unfilled path $\mathcal{S}_{1}$, with central rectangle $\mathcal{R}_{2}$.
	\end{itemize}
\textbf{Step 2:} \ 
	\begin{itemize}
	\item Without loss of generality, assume that the central rectangle ${\mathcal{R}}_{2}$ 
	in $\hat{\mathcal{S}}_{1}\coloneqq \mathcal{S}_{1}$ has a non-integer edge ratio between 
	$1$ and $3$. Otherwise we can use Lemma~\ref{lem:intcase} or Lemma~\ref{lem:expandrec}.
	\item Fill $\mathcal{R}_{2}$ with $\hat{\mathcal{A}}_{2}^{\ast}$, $\hat{\mathcal{B}}_{2}^{\ast}$, 
	$\hat{\mathcal{K}}_{2}^{\ast}$, and $\hat{\mathcal{L}}_{2}^{\ast}$ as described above, leaving 
	the unfilled path $\hat{\mathcal{S}}_{2}$. 
	\item Fill portions of $\mathcal{S}_{1}$ using the sequence defined in the proof of Lemma 
	\ref{lem:fillpath} to create $\mathcal{A}_{2}^{\ast}$, $\mathcal{B}_{2}^{\ast}$, 
	$\mathcal{K}_{2}^{\ast}$, and $\mathcal{L}_{2}^{\ast}$, leaving copies of $\hat{\mathcal{S}}_{2}$. 
	The big dots in Figure~\ref{fig:RhomA2B2} are at the ends of the aforementioned copies of 
	$\hat{\mathcal{S}}_{2}$, thus the union $\mathcal{S}_{2}$ of these copies is a connected 
	unfilled path.
	\end{itemize}
\textbf{Step $\boldsymbol{n}$, $n \geq 3$:} This is treated analogously as Step 2. 
	\begin{itemize}
	\item Without loss of generality, assume that the central rectangle $\mathcal{R}_{n}$ in 
	$\hat{\mathcal{S}}_{n-1}$ has a non-integer edge ratio between $1$ and $3$. Otherwise we 
	can use Lemma~\ref{lem:intcase} or Lemma~\ref{lem:expandrec}.
	\item Fill $\mathcal{R}_{n}$ with $\hat{\mathcal{A}}_{n}^{\ast}$, $\hat{\mathcal{B}}_{n}^{\ast}$, 
	$\hat{\mathcal{K}}_{n}^{\ast}$, and $\hat{\mathcal{L}}_{n}^{\ast}$ as described above, leaving the 
	unfilled path $\hat{\mathcal{S}}_{n}$.
	\item Fill portions of $\mathcal{S}_{n-1}$ using the corresponding sequences of those in the 
	proof of Lemma~\ref{lem:fillpath}, filling out portions of, in order, $\hat{\mathcal{S}}_{n-1}$, 
	$\hat{\mathcal{S}}_{n-2}$, $\ldots$, $\hat{\mathcal{S}}_{2}$, and $\hat{\mathcal{S}}_{1}\coloneqq 
	\mathcal{S}_{1}$, to create $\mathcal{A}_{n}^{\ast}$, $\mathcal{B}_{n}^{\ast}$, $\mathcal{K}_{n}^{\ast}$, 
	and $\mathcal{L}_{n}^{\ast}$. Arguing as before, this leaves an empty path $\mathcal{S}_{n}$ which 
	is composed of copies of $\hat{\mathcal{S}}_{n}$.
	\end{itemize}
 
If $\frac{a}{b}$ is irrational, it is clear that the procedure has infinitely many steps.
Note that if $\frac{a}{b} = \frac{a_{1}}{b_{1}}$  is rational, then the ratios 
$\frac{b_{1}}{v_{1}} = \frac{b_{1}}{|2b_{1}-a_{1}|}$ and 
$\frac{b_{2}}{v_{1}} = \frac{b_{1}}{v_{1}}-q_{1}$ are also rational. In the next step,
we construct arrays for a rectangle with edge ratio $\frac{a_{2}}{b_{2}}$ which
differs from $\frac{v_{1}}{b_{2}}$ by an integer. It follows that the edge
ratio of a rectangle in each succeeding step is rational. Furthermore, the divisor $v_{n}$ 
used at step $n$ is smaller than the remainder $b_{n}$ of the preceding step, so that
the sequence of divisors terminates faster than that of 
the Euclidean algorithm applied to $a$ and $b$. Thus, there are only finitely many steps
in this case.

It is clear from the construction that condition D1 is satisfied. Note also that at each step, 
the squares used to pack the rectangle $\mathcal{R}$ overlap with squares 
used in previous steps only along sets of measure zero. Thus, $\displaystyle\bigcupdot_{i}\mathcal{A}_{i}^{\ast}$ 
and $\displaystyle\bigcupdot_{i}\mathcal{B}_{i}^{\ast}$ are disjoint unions. Furthermore, because 
the width of the unfilled path approaches zero, indeed, condition D2 is satisfied. In fact, if 
$\frac{a}{b}$ is rational, then the width and area of this path becomes zero after a finite number of steps. 
This implies that the unions in condition D2 are finite unions, and it follows that 
$\mathcal{K} \cup \tau(\mathcal{K})$ and $\mathcal{L} \cup \rho(\mathcal{L})$
are disjoint unions up to sets of measure zero. It remains to show that this last condition also holds 
if $\frac{a}{b}$ is irrational. 

Note that for each $t \in \N$, the path that remains unfilled after step $t$ consists of 
$\cl\left(\mathcal{A} \setminus \displaystyle\bigcupdot_{i=1}^{t}\mathcal{A}_{i}^{\ast}\right)$
and $\cl\left(\mathcal{B} \setminus \displaystyle\bigcupdot_{i=1}^{t}\mathcal{B}_{i}^{\ast}\right)$.
Moreover,  $\mathcal{A}_{t+1}^{*} \subseteq 
\cl\left(\mathcal{A} \setminus \displaystyle\bigcupdot_{i=1}^{t}\mathcal{A}_{i}^{\ast}\right)$
and $\mathcal{B}_{t+1}^{*} \subseteq 
\cl\left(\mathcal{B} \setminus \displaystyle\bigcupdot_{i=1}^{t}\mathcal{B}_{i}^{\ast}\right)$. Thus,
$\mathcal{K} = \cl\left(\displaystyle\bigcupdot_{i}\mathcal{K}_{i}^{\ast}\right) = 
\displaystyle\bigcupdot_{i=1}^{t} \mathcal{K}_{i}^{\ast} \cupdot 
\cl\left(\displaystyle\bigcupdot_{i=t+1}^{\infty}\mathcal{K}_{i}^{\ast}\right)$ and
$\mathcal{L} = \cl\left(\displaystyle\bigcupdot_{i}\mathcal{L}_{i}^{\ast}\right) = 
\displaystyle\bigcupdot_{i=1}^{t} \mathcal{L}_{i}^{\ast} \cupdot 
\cl\left(\displaystyle\bigcupdot_{i=t+1}^{\infty}\mathcal{L}_{i}^{\ast}\right)$, where 
$\cl\left(\displaystyle\bigcupdot_{i=t+1}^{\infty}\mathcal{K}_{i}^{\ast}\right) \subseteq 
\cl\left(\mathcal{A} \setminus \displaystyle\bigcupdot_{i=1}^{t}\mathcal{A}_{i}^{\ast}\right)$ and
$\cl\left(\displaystyle\bigcupdot_{i=t+1}^{\infty}\mathcal{L}_{i}^{\ast}\right) \subseteq 
\cl\left(\mathcal{B} \setminus \displaystyle\bigcupdot_{i=1}^{t}\mathcal{B}_{i}^{\ast}\right)$.
Because of condition D1b, it suffices to show that the measure of the unfilled path approaches zero
to justify why $\mathcal{K} \cup \tau(\mathcal{K})$ and $\mathcal{L} \cup \rho(\mathcal{L})$
are disjoint unions up to sets of measure zero.

Consider first the ratio of the measure of $\mathcal{A}_{1}^{\ast} \cup \mathcal{B}_{1}^{\ast}$ to
the measure of $\mathcal{R}$. 

\noindent{Case 1}: Suppose $1<\frac{a_{1}}{b_{1}}<2$. In this case, $q_{1}<\frac{b_{1}}{2b_{1}-a_{1}}<q_{1}+1$, 
and so $\frac{q_{1}+1}{2q_{1}+1} < \frac{b_{1}}{a_{1}} < \frac{q_{1}}{2q_{1}-1}$. 
From the construction, the aforementioned ratio of measures is equal to
\vspace{-7pt} 
$$\begin{aligned}[t]
q_{1}(2q_{1}-1)\cdot\frac{(2b_{1}-a_{1})^2}{a_{1}b_{1}} 
&= q_{1}(2q_{1}-1)\left(\frac{2b_{1}-a_{1}}{b_{1}}\right)^2 \cdot \frac{b_{1}}{a_{1}}
> \frac{q_{1}(2q_{1}-1)}{(q_{1}+1)(2q_{1}+1)} > \frac{1}{6}.
\end{aligned}$$ 

\noindent{Case 2}: Suppose $2<\frac{a_{1}}{b_{1}}<3$. In this case, $q_{1}<\frac{b_{1}}{a_{1}-2b_{1}}<q_{1}+1$, 
and so $\frac{q_{1}}{2q_{1}+1} < \frac{b_{1}}{a_{1}} < \frac{q_{1}+1}{2q_{1}+3}$. 
From the construction, the aforementioned
ratio of measures is equal to 
\vspace{-7pt}
$$\begin{aligned}[t]
q_{1}(2q_{1}+1)\cdot\frac{(a_{1}-2b_{1})^2}{a_{1}b_{1}} 
&= q_{1}(2q_{1}+1)\left(\frac{a_{1}-2b_{1}}{b_{1}}\right)^2\cdot\frac{b_{1}}{a_{1}}
> \left(\frac{q_{1}}{q_{1}+1}\right)^2 > \frac{1}{4}.
\end{aligned}$$ 

We note that in the procedure outlined in the proof of Lemma~\ref{lem:fillpath}, each of the 
last two rectangles overlap with some previous rectangle except on a square region congruent to 
$\hat{\mathcal{A}}$. Thus, we also look for a lower bound for the ratio of the measure of 
$\mathcal{A}_{1}^{\ast}$ to that of $\mathcal{A}$. In either of the two cases above,
this ratio is
\vspace{-7pt}
$$\begin{aligned}[t]
q_{1}^2 \cdot \frac{|2b_{1}-a_{1}|^2}{b_{1}^{2}} &  > \left(\frac{q_{1}}{q_{1}+1}\right)^2 > \frac{1}{4}.
\end{aligned}$$ 

These computations show that at each step of the construction, the ratio of the measure of the portion 
that will remain unfilled is at most $\frac{5}{6}$ of the measure of the unfilled portion in the 
previous step. Thus, the measure of the unfilled portion approaches zero as the construction is 
carried out. It follows that $\mathcal{K}$ and $\mathcal{L}$ satisfy Rh1, Rh2, and Rh3.

Forming $\mathcal{S}$ as in Figure~\ref{fig:RhomStepQ} and 
$\mathcal{E} = P(\Gamma)(\mathcal{K} \cupdot \mathcal{L})$, we have that 
$F = \cl(\mathcal{S} \cup \mathcal{E})$ is a compact fundamental domain for $\Gamma$ 
with $[S(F):P(\Gamma)]=2$. This completes the proof of Theorem~\ref{thm:rho}.
\end{proof}

\begin{figure}[ht]
\begin{center}
\includegraphics[width=0.95\textwidth]{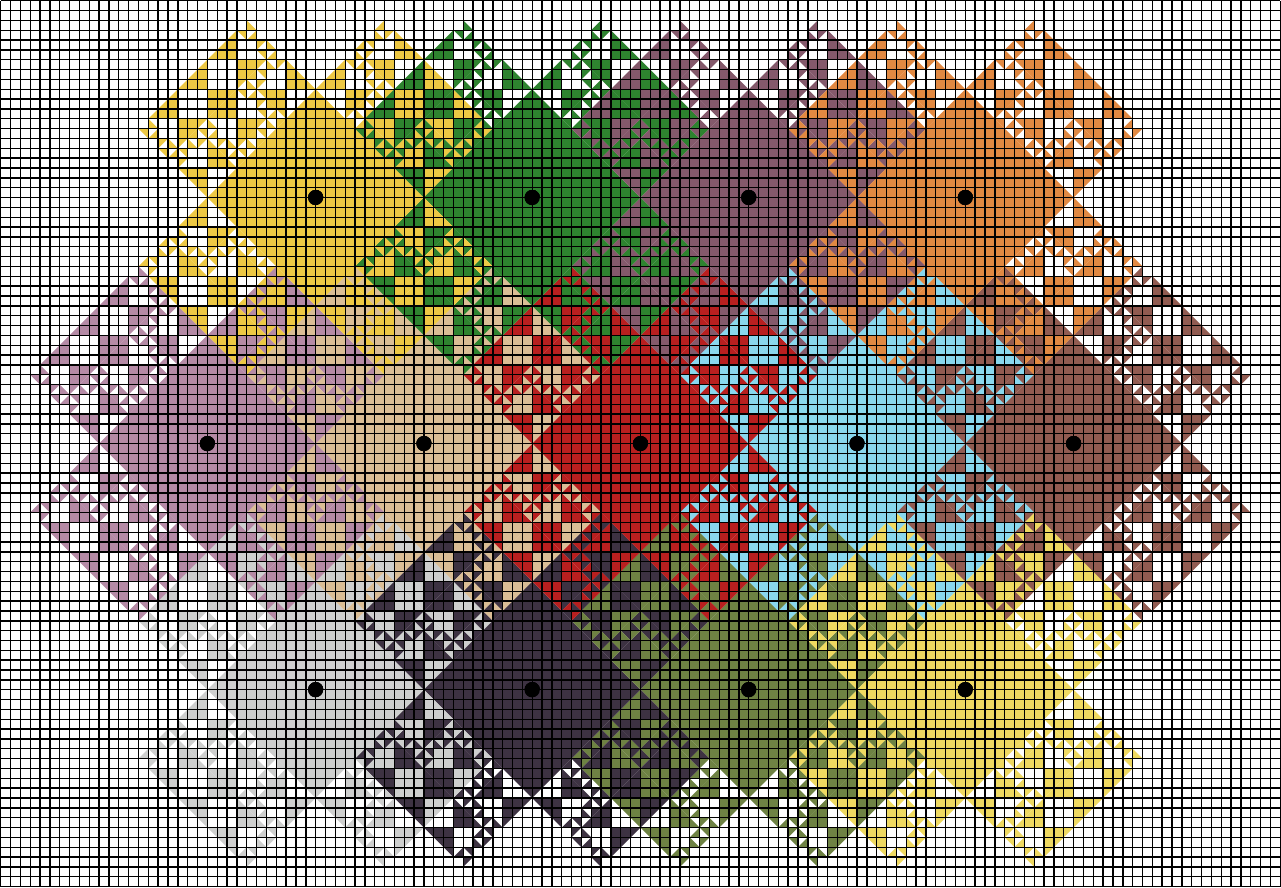}
\caption{Fundamental domain for the rhombic lattice $\Gamma$ with $m=11$, $n=25$,
constructed using the procedure, together with some of its $\Gamma$-translates.}
\label{fig:RhomFDEx}
\end{center}
\end{figure}

In Figure~\ref{fig:RhomFDEx}, we show a portion of the tiling by the fundamental domain 
for the rhombic lattice $\Gamma$ with basis $\left\{\vect{22}{0},\vect{11}{25}\right\}$ 
generated by the procedure described.

%%%%%%%%%%%%%%%%%%%%%%%%%%%%%%%%%%%%%%%%%%%%%%%%%%%%%%%
\subsection*{Acknowledgements}

The authors would like to thank the anonymous referee for valuable remarks that improved 
the quality of the paper. J.R.C.G.~Damasco is grateful to the University of the Philippines 
System for financial support through its Faculty, REPS, and Administrative Staff Development Program. 
D.~Frettl\"oh is grateful to the Research Center for Mathematical Modelling, Bielefeld 
University, and the University of the Philippines System for financial support.

%%%%%%%%%%%%%%%%%%%%%%%%%%%%%%%%%%%%%%%%%%%%%%%%%%%%%%%
% You do not have to use the same format for your references, but 
%    include everything in this file.  Don't use natbib please.
% If you use BibTeX to create a bibliography, copy the .bbl file into here.


\begin{thebibliography}{99}

\bibitem{bks}
M.\ Baake, R.\ Klitzing, M.\ Schlottmann.
Fractally shaped acceptance domains of quasiperiodic square-triangle 
tilings with dodecagonal symmetry.
{\em Physica A}, 191:554-558, 1992.

\bibitem{coc}
E.\ Cockayne. 
Nonconnected atomic surfaces for quasicrystalline sphere packings.
{\em Phys.\ Rev.\ B}, 49:5896-5910, 1994.

\bibitem{df}
J. Damasco, D. Frettl\"oh, M. Loquias. 
Highly symmetric fundamental domains for lattices in $\R^2$ and $\R^3$.
{\tt arXiv:1305.1798 [math.CO]}.

\bibitem{dehn} 
M. Dehn. \"Uber Zerlegung von Rechtecken in Rechtecke. {\em Math. Annalen},
57:314-332, 1903.

\bibitem{elser}
V.\ Elser.
Exceptionally symmetric fundamental domains for the root lattices in 2D. 
In {\em Aperiodic Order}, Tagungsbericht Oberwolfach, No.~20, Mathematisches 
Forschungsinstitut Oberwolfach, Germany, 2001. 6.

\bibitem{elserweb}
V.\ Elser. Fractal Fun. \texttt{http://uuuuuu.lassp.cornell.edu/gallery/fractal\_fun}, 
accessed 05 January 2018.

\end{thebibliography}
\end{document}